\documentclass{amsart}
\usepackage{amssymb, mathrsfs, bbm}
\usepackage[T1]{fontenc}
\usepackage[sc, osf]{mathpazo}
\selectfont
\usepackage[protrusion=true, expansion=true]{microtype}

\usepackage[all, arc]{xy}
\SelectTips{eu}{}
\entrymodifiers={+!!<0pt,\fontdimen22\textfont2>}

\usepackage[pdftex, colorlinks=true, linkcolor=blue, citecolor=blue, linktocpage]{hyperref}

\makeatletter
\def\tagform@#1{\maketag@@@{\bfseries(\ignorespaces#1\unskip\@@italiccorr)}}
\renewcommand{\eqref}[1]{\textup{{\normalfont(\ref{#1}}\normalfont)}}
\makeatother

\numberwithin{equation}{subsection}
\theoremstyle{plain}
\newtheorem{thm}[equation]{Theorem}

\theoremstyle{definition}
\newtheorem{example}[equation]{Example}

\theoremstyle{remark}

\theoremstyle{definition}

\def\11{\mathbf{1}}

\def\AA{\mathbf{A}}

\def\GG{\mathbf{G}}

\def\PP{\mathbf{P}}

\def\ZZ{\mathbf{Z}}

\def\dim{\mathrm{dim}}

\def\res{\mathrm{res}}

\newcommand{\mapright}[1]{\xrightarrow{#1}}

\title{Motivic splitting principle}
\author{Rahbar Virk}
\address{The Appalachians}
\begin{document}
\maketitle
\renewcommand{\thesubsection}{\arabic{subsection}}
\subsection{Introduction}
Fix a field $k$ and a commutative ring $\Lambda$ (with $1$). Assume either that $k$ is perfect and admits resolution of singularities, or that $k$ is arbitrary and the exponential characteristic\footnote{
If $\mathrm{char}(k)=0$, then the exponential characteristic of $k$ is $1$. If $\mathrm{char}(k)> 0$, then the exponential characteristic of $k$ is $\mathrm{char}(k)$.}
of $k$ is invertible in $\Lambda$. Write $DM(k; \Lambda)$ for the triangulated category of motives with $\Lambda$-coefficients (see \S\ref{s:motives}).

We write `scheme' in lieu of `separated scheme of finite type over $k$'. For a group scheme $G$, a \emph{$G$-torsor} will mean a faithfully flat $G$-invariant morphism $X\to Y$, such that the canonical morphism $G\times X \to X\times_Y X$ is an isomorphism. In this situation, we say that the \emph{quotient} of the $G$-action on $X$ exists, and we set $X/G = Y$. The \emph{trivial torsor} over a scheme $X$ is the projection $G\times X \to X$, with $G$ acting on $G\times X$ via multiplication on $G$.
A \emph{reductive group} will mean a smooth affine group scheme $G$ such that every smooth connected unipotent subgroup of $G\times_k \bar k$ is trivial, where $\bar k$ is the algebraic closure of $k$.

\begin{thm}\label{splittingThm}Let $G$ be a connected split reductive group, and let $B\subset G$ be a Borel subgroup.
Let $X$ be a scheme with $G$-action.
Let $t(G)$ denote the torsion index of $G$. If $t(G)$ is invertible in $\Lambda$, then there is an isomorphism
\[ M^c(X_{hB}) \simeq M^c(X_{hG}) \otimes M^c(G/B)\]
in $DM(k;\Lambda)$ that commutes with smooth pullbacks and localization triangles.
\end{thm}
Here $M^c(X_{hG})$, $M^c(X_{hB})$ denote the $G$-equivariant and $B$-equivariant (Borel-Moore) motives, respectively, of $X$ (in the sense of \cite{T}; see \S\ref{s:equi}).

The proof of Theorem \ref{splittingThm} is essentially the same as that of the analogous statement for the cohomology of topological spaces (for instance, compare with \cite{May} and \cite[Theorem 16.1]{TB}). No claim to originality is being made. 
%

\begin{example}If $X=\mathrm{Spec}(k)$, then Theorem \ref{splittingThm} states
\[ M^c(BB) \simeq M^c(BG)\otimes M^c(G/B), \]
where $M^c(BB)$ and $M^c(BG)$ are the motives of the classifying spaces of $B$ and $G$, respectively (see \cite{T} and \S\ref{s:equi}). This is a motivic analogue of the usual splitting principle (working with groups over the complex numbers say):
\[ H^*(BB) \simeq H^*(BG) \otimes H^*(G/B), \]
where cohomology is with coefficients in a ring $\Lambda$ in which $p$ is invertible for all primes $p$ such that $H^*(G; \ZZ)$ has $p$-torsion.
\end{example}

\subsection{Motives}\label{s:motives}
The category $DM(k; \Lambda)$ is the monoidal triangulated category of Nisnevich motivic spectra over $k$, obtained by applying $\AA^1$-localization and $\PP^1$-stabilization to the derived category of (unbounded) complexes of Nisnevich sheaves with transfer \cite[Definition 11.1.1]{CD}.\footnote{
The category $DM^{eff}_-(k;\Lambda)$ of \cite{V} embeds into $DM(k;\Lambda)$ as a full and faithful subcategory \cite[Example 11.1.3]{CD}. 
Utilizing $DM(k;\Lambda)$ (instead of $DM^{eff}_-(k;\Lambda)$) is dictated by the need for arbitrary direct sums and compact generation \cite[Corollary 3.5.5]{V}. Both of these properties are required to make sense of equivariant motives (see \cite{T} and \S\ref{s:equi}).}
The primary references for the properties of $DM(k;\Lambda)$ that we use are \cite{V}, \cite{CD} and \cite{CD1}.
The assumption that $k$ is perfect and admits resolution of singularities stems from the treatment in \cite{V}. The alternate assumption that $k$ is arbitrary, but the exponential characteristic is invertible in $\Lambda$, stems from \cite[Proposition 8.1]{CD1} and \cite{K}. These articles, amongst other things, extend the constructions of \cite{V} under this alternate hypothesis.

There is a covariant functor $X\mapsto M^c(X)$ from the category of schemes and proper morphisms to $DM(k;\Lambda)$ (see \cite[\S2.2]{V}; alternatively, in the notation of \cite{CD}, we have $M^c(X) = a_*a^!\Lambda$, where $a\colon X \to \mathrm{Spec}(k)$ is the structure morphism). The functor $M^c(X)$ behaves like a Borel-Moore homology theory in the following sense. If $f\colon X\to Y$ is a smooth morphism with fibres of dimension $r$ (the morphism $f$ is allowed to have some fibres empty), then there is a map (compatible with composition of morphisms):
\[ f^*\colon M^c(Y)(r)[2r] \to M^c(X),\]
where $(j)$ denotes tensoring with the $j$-th Tate twist (see \cite[(11.1.2.2]{CD}) and $[i]$ denotes the $i$-th shift functor (available in any triangulated category).
If $f$ is a vector bundle, then $f^*$ is an isomorphism. 

Let $i\colon Z \hookrightarrow X$ be a closed immersion, and let $j\colon U\hookrightarrow X$ be the open immersion of the complement $U=X-Z$. Then $i_*$ and $j^*$ 
fit into a canonical
distinguished triangle \cite[\S2.2]{V}, \cite[Corollary 5.9, Theorem 5.11]{CD1}, the \emph{localization triangle},
\[ M^c(Z) \mapright{i_*} M^c(X) \mapright{j^*} M^c(U) \mapright{\partial_i}\]
The category $DM(k;\Lambda)$ is a symmetric monoidal triangulated category, and
\[ M^c(X \times Y) = M^c(X) \otimes M^c(Y).\]
The motive $M^c(\mathrm{Spec}(k))$ is the unit object. 
Although notationally abusive, it is convenient to set
\[ \Lambda = M^c(\mathrm{Spec}(k)).\]
Let $H^i_M(X; \Lambda(j))$ denote the motivic cohomology groups of $X$, as defined in \cite[\S11.2]{CD}. These are contravariant functors from the category of schemes to $\Lambda$-modules. By \cite[Example 11.2.3]{CD}, if $X$ is smooth and equidimensional, then motivic cohomology determines the Chow ring $CH^*(X)$ of $X$:
\[ H^{2j}_M(X; \Lambda(j)) = CH^j(X) \otimes_{\ZZ} \Lambda.\]
For an arbitrary scheme $X$, each $e\in H^i(X;\Lambda(j))$ determines a canonical map
\[ e\cap\colon M^c(X)(-j)[-i] \to M^c(X).\]
\begin{example}
Let $f\colon L\to X$ be a line bundle. Write $i\colon X \hookrightarrow L$ for the zero section. Let $c_1(L) \in H^2_M(X; \Lambda(1))$ be the first Chern class of $L$ \cite[Definition 11.3.2]{CD}. Then 
$c_1(L)\cap$
is the composition
\[ M^c(X)(-1)[-2] \mapright{i_*}M^c(L)(-1)[-2] \mapright{f^{*-1}} M^c(X).\]
\end{example}

\subsection{Equivariant motives, following B. Totaro}\label{s:equi}
Let $G$ be an affine group scheme. Let $X$ be a scheme with $G$-action. Let
$\cdots \twoheadrightarrow V_2 \twoheadrightarrow V_1$
be a sequence of surjections of $G$-equivariant vector bundles over $X$. Write $n_i$ for the rank of $V_i$. Let $U_i\subset V_i$ be a $G$-stable open subscheme such that $V_{i+1}-U_{i+1}$ is contained in the inverse image of $V_i-U_i$, and such that the quotient $U_i/G$ exists. Assume that the codimension of $V_i-U_i$ goes to infinity with $i$. The \emph{equivariant motive} $M^c(X_{hG})$ is defined to be the homotopy limit of the sequence
\[ \cdots \to M^c(U_2/G)(-n_2)[-2n_2]\to M^c(U_1/G)(-n_1)[-2n_1].\]
The definition of $M^c(X_{hG})$ is independent (up to a not necessarily unique isomorphism) of all the choices involved \cite[Theorem 8.5]{T}.
If $X$ satisfies any of the conditions of \cite[Proposition 23]{EG}, then such vector bundles exist (also see \cite[Remark 1.4]{To1} and \cite[\S8]{T}). Whenever we speak of $M^c(X_{hG})$, we implicitly, and without further comment, assume that such vector bundles exist (including in the statement of Theorem \ref{splittingThm}).

Since $M^c(X_{hG})$ is defined as a (homotopy) limit of ordinary motives $M^c(U_i/G)$, the functorial properties of ordinary motives (pullback, pushforward, localization triangles, etc.) extend to the equivariant setup. By construction, $M^c(X_{hG})$ satisfies equivariant descent: if the quotient of the $G$-action on $X$ exists, then
\[ M^c(X_{hG}) \simeq M^c(X/G).\]
We set
\[ M^c(BG) = M^c(\mathrm{Spec}(k)_{hG}).\]
\begin{example}
Let $V_i$ be the direct sum of $i$-copies of the natural $1$-dimensional representation of $\GG_m$. Let $U_i = V_i - \{0\}$. Then $U_i/\GG_m\simeq\PP^{i-1}$. We have \cite[Lemma 8.7]{T}:
\[ M^c(B\GG_m) \simeq \prod_{i\leq -1}\Lambda(i)[2i].\]
\end{example}

\subsection{Restriction to a subgroup}Let $G$ be an affine group scheme. Let $H\subset G$ be a closed subgroup. If the quotient $X/G$ exists, then the quotient $X/H$ exists. If $G$ is smooth, then we have a pullback
\[ M^c(X/G)(\dim(G/H))[2\dim(G/H)] \to M^c(X/H). \]
This family of pullbacks, one for each such $X$, yields a map, \emph{restriction},
\[ \res^H_G\colon M^c(Y_{hG})(\dim(G/H))[2\dim(G/H)] \to M^c(Y_{hH}),\]
for any scheme $Y$ with $G$-action.
Restriction commutes with smooth pullbacks and localization triangles.

\subsection{Chow ring of a classifying space}
Let $G$ be an affine group scheme. Following \cite{To1}, define the Chow ring $CH^*_G$ of the classifying space of $G$ as follows. Let $V$ be a representation of $G$ over $k$. Let $U\subset V$ be an open subscheme such that the quotient $U/G$ exists, and such that $V-U$ has codimension greater than $i$. Then $CH^i_G = CH^i(U/G)$. This definition is independent of all the choices involved \cite[Theorem 1.1]{To1} and gives a well-defined ring $CH^*_G$.
It follows from the definition that each $e\in CH^i_G$ determines a canonical map
\[ e\cap \colon M^c(X_{hG})(-i)[-2i] \to M^c(X_{hG}),\]
for a scheme $X$ with an action of $G$. 

By faithfully flat descent, each representation of $G$ over $k$ determines a vector bundle over the schemes $U/G$ used to define $CH^*_G$. Consequently, each such representation has Chern classes in $CH^*_G$.
\begin{example}
Let $T$ be a split torus.
Let $\chi$ be a character of $T$, with first Chern class $c_1(\chi) \in CH^1_{T}$. Let $X$ be a scheme with $T$-action. Then $\chi$ determines an equivariant line bundle $L_{\chi}\to X$. If the quotient $X/T$ exists, then $L_{\chi}$ descends to a line bundle $\tilde L_{\chi}\to X/T$. In this situation, the map
\[ c_1(\chi)\cap\colon M^c(X_{hT})(-1)[-2] \to M^c(X_{hT}) \]
is the composition
\[ M^c(X_{hT})(-1)[-2] \mapright{\sim}M^c(X/T)(-1)[-2] \mapright{c_1(\tilde L_{\chi})\cap} M^c(X/T) \mapright{\sim}M^c(X_{hT}),\]
where the first and last isomorphisms are taken to be inverse to each other.
\end{example}

\subsection{The torsion index}\label{s:tor}
Let $G$ be a connected split reductive group over $k$.
Let $B\subset G$ be a Borel subgroup. The \emph{torsion index} of $G$ is the smallest integer $t(G)\in \ZZ_{>0}$ such that the image of the map $CH^*_B\to CH^*(G/B)$ contains $t(G)\cdot CH^{\dim(G/B)}(G/B)$. The natural map,
\[ CH^*_B \otimes_{\ZZ} \ZZ[t(G)^{-1}] \to CH^*(G/B) \otimes_{\ZZ} \ZZ[t(G)^{-1}],\]
is surjective.
According to \cite[Th\'eor\`eme 2]{G}, for any $G$-torsor $X\to Y$, there is a non-empty open subscheme $U\subset Y$ along with a finite \'etale morphism $V\to U$ of degree invertible in $\ZZ[t(G)^{-1}]$, such that $X$ is trivial over $V$.\footnote{
The point is that, if one stays away from primes that divide $t(G)$, then all the challenges of `\'etale descent' for equivariant Chow groups disappear.
For further information on the torsion index, \cite{ToSpin} is highly recommended.}
\begin{example}
The group $GL_n$ has torsion index $1$. 
\end{example}
\subsection{Proof of Theorem \ref{splittingThm}}(Compare with \cite{May} and the proof of \cite[Theorem 16.1]{TB}).
Pick elements $e_1, \ldots, e_{n}\in CH^*(BB) \otimes_{\ZZ} \Lambda$, of homogeneous degree, that restrict to a basis of $CH^*(G/B)\otimes_{\ZZ}\Lambda$. 
Write $d_i$ for the degree of $e_i$. Set $d=\dim(G/B)$.
For each $e_i$, consider the composition
\[ M^c(X_{hG})(d-d_i)[2(d-d_i)] \mapright{\res^B_G} M^c(X_{hB})(-d_i)[-2d_i] \mapright{e_i\cap } M^c(X_{hB}).\]
Summing these, we obtain a map
\[ \bigoplus_{i}M^c(X_{hG})(d-d_i)[2(d-d_i)] \to M^c(X_{hB}). \]
By the Bruhat decomposition, this may be rewritten as a map
\[ \theta\colon M^c(X_{hG}) \otimes M^c(G/B) \to M^c(X_{hB}).\]
The map $\theta$ commutes with smooth pullbacks and localization triangles. 
We will show $\theta$ is an isomorphism.
It suffices to demonstrate this under the assumption that the quotient $X/G$ exists. Via the isomorphisms $M^c(X_{hG}) \simeq M^c(X/G)$ and $M^c(X_{hB}) \simeq M^c(X/B)$, the map $\theta$ yields a map
\[ \theta_X\colon M^c(X/G)\otimes M^c(G/B) \to M^c(X/B). \]
If $X\to X/G$ is the trivial $G$-torsor, then $\theta_X$ is manifestly an isomorphism.
In general, there exists a non-empty open subscheme $U\subset X/G$, along with a finite \'etale morphism $f\colon V\to U$ of degree invertible in $\ZZ[t(G)^{-1}]$, 
such that $X$ pulled back to $V$ is the trivial $G$-torsor (see \S\ref{s:tor}). 
The map $f_*f^*\colon M^c(U)\to M^c(U)$ is the degree of $f$ times the identity (as follows from \cite[A.5 (6)]{CD} and \cite[Proposition 11.2.5]{CD}). 
Consequently, $\theta_U$ is an isomorphism. Now let $Z=X-U$ be the closed complement (with reduced scheme structure say). Then, by virtue of the localization triangle, it suffices to show $\theta_Z$ is an isomorphism. This follows from an induction on dimension (the base case has been dealt with by the above considerations). 

\subsection{A complement}
The Chow ring $CH^*_B$ acts on $M^c(X_{hB})$. Under the isomorphism
\[ M^c(X_{hB}) \simeq M^c(X_{hG}) \otimes M^c(G/B), \]
this action on the right hand side is the action of $CH^*_B$ on $M^c(G/B)$.

\end{document}